\documentclass[12pt]{tac}

\usepackage{amsfonts,amsmath}
\usepackage{paralist}
\usepackage{bbm}
\usepackage[all]{xy}

\title{A Presentation of the Category of Stochastic Matrices}
\author{Tobias Fritz}
\thanks{I am grateful to the Max Planck Institute for providing an excellent research environment and finanical support.}
\address{Max Planck Institute for Mathematics\\ Vivatsgasse 7, 53111 Bonn, Germany}
\copyrightyear{2009}
\eaddress{fritz@mpim-bonn.mpg.de}

\newcommand{\N}{\mathbb{N}}
\newcommand{\R}{\mathbb{R}}
\newcommand{\beq}{\begin{displaymath}}
\newcommand{\eeq}{\end{displaymath}}
\newcommand{\beqn}{\begin{equation}}
\newcommand{\eeqn}{\end{equation}}
\newcommand{\beqa}{\begin{eqnarray*}}
\newcommand{\eeqa}{\end{eqnarray*}}
\newcommand{\ra}{\rightarrow}
\newcommand{\lra}{\longrightarrow}
\newcommand{\id}{\mathrm{id}}
\newcommand{\re}[1]{~(\ref{#1})}
\newtheorem{prop}{Proposition}[section]
\newtheorem{lem}[prop]{Lemma}
\newtheorem{thm}[prop]{Theorem}
\theoremstyle{definition}
\newtheorem{defn}[prop]{Definition}
\theoremstyle{remark} 

\begin{document}

\maketitle
\begin{abstract}
This note gives generators and relations for the strict monoidal category of probabilistic maps on finite cardinals (i.e., stochastic matrices).
\end{abstract}

\setcounter{section}{-1}
\section{Notation}

$\mathtt{FinMap}$ is the category of finite cardinals with ordinary functions as morphisms. The notation $[n]$ is shorthand for the $n$-element set $\{1,\ldots,n\}$ and is identified with the corresponding cardinal. $\mathbbm{1}_n$ denotes the unit matrix of size $n\times n$. The acronym ``i.a.'' stands for ``induction assumption''.

\section{Introduction}

Algebraic structures like groups, rings or lattices can be defined via their universal instances, the so-called \emph{Lawvere theories}. Recall that a Lawvere theory is a category $\mathtt{LT}$ with finite products together with a product-preserving functor $\mathtt{FinMap}^\mathrm{op}\ra\mathtt{LT}$ which is bijective on objects. 

Usually, one defines an algebraic structure in terms of a family of operations of specific arity. Then this family of operations together with the structure-defining equations between them forms a presentation of the corresponding Lawvere theory $\mathtt{LT}$. However, in other cases it may happen that we have $\mathtt{LT}$ defined directly as a category, and we want to recover a family of operations together with a family of equations between these, such that this data defines the same algebraic structure as $\mathtt{LT}$ does. This is equivalent to determining a presentation of $\mathtt{LT}$, and this is what will be done here for the particular case $\mathtt{LT}^\mathrm{op}=\mathtt{FinStoMap}$, where $\mathtt{FinStoMap}$ is the category of ``probabilistic maps'' on finite cardinals (see below). For the reason of calculational simplicity, the given presentation is a presentation of $\mathtt{FinStoMap}$ as a strict monoidal category with respect to the coproduct, and not a presentation of $\mathtt{FinStoMap}$ as a category with finite coproducts.

For other examples of presentations of Lawvere theories as strict monoidal categories, see~\cite{Laf}. That article in particular contains a presentation of $\mathtt{FinMap}$ given by the first three generators and the first five relations of definition~\ref{pres} together with the equations\re{edelisid} and\re{psp}. Although the present article is self-contained, knowledge of~\cite{Laf} will help in understanding the proofs presented here.

The present results are applied in~\cite{Fri} to the study of \emph{convex spaces}, which are an abstract version of convex subsets of vector spaces. A priori, a convex space is a model of $\mathtt{FinStoMap}^\mathrm{op}$. Theorem~\ref{FinStoMapPres} however facilitates a description of convex space structure in terms of a family of binary convex combination operations satisfying certain compatibility conditions.

\section{The category of stochastic matrices}

In this article, the term \emph{stochastic matrix} means \emph{column-stochastic matrix over $\R$}, i.e. a matrix with nonnegative real entries such that each column sums to $1$. The product of two stochastic matrices is again a stochastic matrix. One way to think of a stochastic matrix $A$ of size $n\times m$ is as a probabilistic map $[m]\ra[n]$, meaning that it assigns to every $j\in[m]$ a probability distribution on $[n]$, and these assignments are probabilistically independent. It is useful to visualize this process as a braid-like diagram\\
\beqn
\label{blackbox}
\vcenter
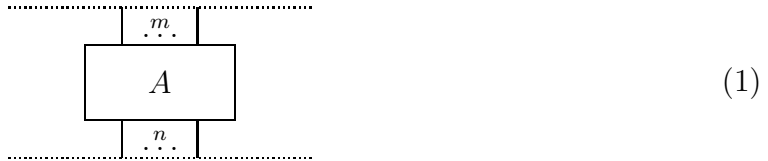
{\xy
(0,0);(40,0) **@{.},
(0,20);(40,20) **@{.},
(10,5);(30,5) **@{-},
(10,15);(30,15) **@{-},
(20,10)*={A};
(10,5);(10,15) **@{-},
(30,5);(30,15) **@{-},
(15,0);(15,5) **@{-},
(25,0);(25,5) **@{-},
(20,2.5)*={\stackrel{n}{\ldots}};
(15,15);(15,20) **@{-},
(25,15);(25,20) **@{-},
(20,17.5)*={\stackrel{m}{\ldots}};
\endxy}
\eeqn\\
with $m$ input strands, representing the elements of $[m]$, and $n$ output strands, representing the elements of $[n]$, and a picture of the strands crossing, coalescing, and newly emerging, here drawn as a ``black box'' $A$. In case of a deterministic map $[m]\ra[n]$, each of the $m$ input strands would get mapped to a unique output strand. However, now in the case of probabilistic maps, an input strand may also branch into several output strands, where each branch carries a certain fraction of the input strand. 

As a degenerate case, we stipulate that there exists exactly one stochastic matrix of size $n\times 0$ for each $n$, corresponding to the unique function $[0]=\emptyset\ra[n]$.

\begin{defn}[The finitary stochastic map category $\mathtt{FinStoMap}$]
\beqa
\mathrm{Obj}(\mathtt{FinStoMap})&\equiv&\N_0\textrm{ (finite cardinals)}\\
\mathtt{FinStoMap}(m,n)&\equiv&\textrm{stochastic matrices of size $n\times m$}
\eeqa
Composition is defined by matrix multiplication.
\end{defn}

It is clear that this satisfies the axioms of a category, as matrix multiplication is associative and the unit matrices act as identity morphisms. In the diagram picture, composition is represented by vertical juxtaposition of the diagrams.

As an equivalent definition, one might take the morphisms in $\mathtt{FinStoMap}$ to be the conditional probability distributions on $[n]$ dependent on a distribution on $[m]$. Composition is then given by the Chapman-Kolmogorov equation. A third formulation could be as the category of communication channels on finite alphabets with concatenation of channels as composition of morphisms. 

The goal of this article is to find a different and purely algebraic description of $\mathtt{FinStoMap}$ in terms of generators and relations with respect to the strict monoidal structure given by the coproduct. This is related to but more elaborate than giving a presentation of
\begin{compactitem}
\item a symmetric group $S_n$ (see for example~\cite[6.2]{Cox})
\item the category $\mathtt{FinMap}$ (see~\cite{Laf} for a precise statement and proof)
\end{compactitem}
Simpler variants of the statements and proofs given here would also apply to yield the cited standard solutions to these two problems.

\begin{lem}
$\mathtt{FinStoMap}$ has all finite coproducts.
\end{lem}

\begin{proof}
$0\in\mathrm{Obj}(\mathtt{FinStoMap})$ clearly is an initial object, thereby defining the empty coproduct. Now for binary coproducts of two objects $[n_1]$ and $[n_2]$. The inclusion morphisms are 
\beq
\left(\begin{array}{c}\mathbbm{1}_{n_1}\\0\end{array}\right):[n_1]\ra[n_1+n_2],\qquad\left(\begin{array}{c}0\\\mathbbm{1}_{n_2}\end{array}\right):[n_2]\ra[n_1+n_2].
\eeq
They satisfy the universal property
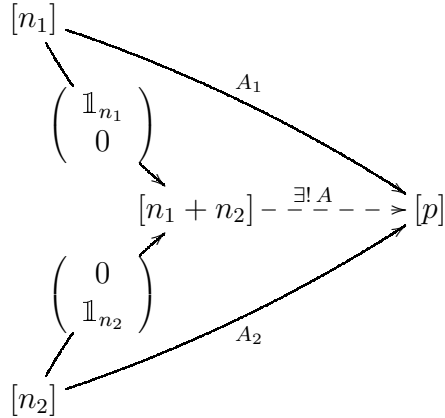
\beq
\xymatrix{[n_1]\ar@/^/[rrrdd]^{A_1}\ar@/_/[rdd]|{\left(\begin{array}{c}\mathbbm{1}_{n_1}\\0\end{array}\right)} \\\\
& [n_1+n_2]\ar@{-->}[rr]^{\exists!\,A} && [p] \\\\
[n_2]\ar@/_/[rrruu]_{A_2}\ar@/^/[ruu]|{\left(\begin{array}{c}0\\\mathbbm{1}_{n_2}\end{array}\right)}}
\eeq
since commutativity of this diagram is equivalent to $A=\left(\begin{array}{cc}A_1&A_2\end{array}\right)$. This $A$ is clearly a stochastic matrix provided that both $A_1$ and $A_2$ are.
\end{proof}

In the following, $\mathtt{FinStoMap}$ will be regarded as a strict monoidal category with respect to the coproduct. Then the monoidal product of two stochastic matrices $A_1:[m_1]\ra[n_1]$ and $A_2:[m_2]\ra[n_2]$ is the block-diagonal matrix
\beq
\left(\begin{array}{cc}A_1&0\\0&A_2\end{array}\right):[m_1+m_2]\ra[n_1+n_2].
\eeq
In the ``black box'' picture, this product is represented by horizontal juxtaposition of diagrams. Note that when $m_1=0$, the resulting matrix is just $A_2$, together with an additional collection of rows only containing zeros. Similarly when $m_2=0$.

\section{$\mathtt{FinStoMap}$ by generators and relations}

What follows now is the definition of a strict monoidal category $\mathtt{FinStoMap'}$ in terms of generators and relations. In the following definition, domain and codomain of each generator are indicated by the number of input strands and output strands, respectively, of each diagrammatic representation.

\begin{defn}
\label{pres}
$\mathtt{FinStoMap'}$ is the strict monoidal category generated by one object $[1]$ with tensor powers $[n]=[1]^{\otimes n}$ together with the family of morphisms
\beqa
\vcenter{\xy
(-10,10)*+{\partial=};
(0,0)="LD";(40,0)="RD";(0,20)="LU";(40,20)="RU";
"LD";"RD" **@{.},
"LU";"RU" **@{.},
(20,0)="O1";(20,10)*={\times}="C1";
"C1";"O1" **@{-},
(50,10)*={e=};
(60,0)="LD";(100,0)="RD";(60,20)="LU";(100,20)="RU";
"LD";"RD" **@{.},
"LU";"RU" **@{.},
(70,20)="I1";(90,20)="I2";(80,10)="C1";(80,0)="O1";
"C1";"O1" **@{-},
"I1";"C1" **\crv{(70,15)&(80,15)},
"I2";"C1" **\crv{(90,15)&(80,15)},
\endxy}
\\\\\\
\vcenter{\xy
(-10,10)*={s=};
(0,0)="LD";(40,0)="RD";(0,20)="LU";(40,20)="RU";
"LD";"RD" **@{.},
"LU";"RU" **@{.},
(10,0)="O1";(30,0)="O2";(10,20)="I1";(30,20)="I2";
"O1";"I2" **\crv{(10,10)&(30,10)},
"O2";"I1" **\crv{(30,10)&(10,10)},
(50,10)*={c_\lambda=};
(60,0)="LD";(100,0)="RD";(60,20)="LU";(100,20)="RU";
"LD";"RD" **@{.},
"LU";"RU" **@{.},
(80,20)="I1";(70,0)="O1";(90,0)="O2";(80,10)*+[F]{\lambda}="C1";
"O1";"C1" **\crv{(70,5)&(80,5)},
"O2";"C1" **\crv{(90,5)&(80,5)},
"I1";"C1" **@{-},
(100,10)*={\forall\lambda\in[0,1]};
\endxy}\eeqa
subject to the relations
\beqn
\label{eass}
\vcenter{\xy
(20,25)*={e(e\otimes\mathrm{id}_{[1]})};
(50,25)*={=};
(80,25)*={e(\mathrm{id}_{[1]}\otimes e)\::};
(0,0);(40,0) **@{.},
(0,10);(40,10) **@{.},
(0,20);(40,20) **@{.},
(10,20);(15,15) **\crv{(10,17.5)&(15,17.5)},
(20,20);(15,15) **\crv{(20,17.5)&(15,17.5)},
(15,15);(15,10) **@{-},
(30,20);(30,10) **@{-},
(22.5,0);(22.5,5) **@{-},
(22.5,5);(15,10) **\crv{(22.5,7.5)&(15,7.5)},
(22.5,5);(30,10) **\crv{(22.5,7.5)&(30,7.5)},
(50,10)*={=};
(60,0);(100,0) **@{.},
(60,10);(100,10) **@{.},
(60,20);(100,20) **@{.},
(70,20);(70,10) **@{-},
(80,20);(85,15) **\crv{(80,17.5)&(85,17.5)},
(90,20);(85,15) **\crv{(90,17.5)&(85,17.5)},
(85,15);(85,10) **@{-},
(70,10);(77.5,5) **\crv{(70,7.5)&(77.5,7.5)},
(85,10);(77.5,5) **\crv{(85,7.5)&(77.5,7.5)},
(77.5,5);(77.5,0) **@{-}
\endxy}
\eeqn\\
\beqn
\label{ese}
\vcenter{\xy
(20,25)*={es};
(50,25)*={=};
(80,25)*={e\::};
(0,0);(40,0) **@{.},
(0,10);(40,10) **@{.},
(0,20);(40,20) **@{.},
(20,0);(20,5) **@{-},
(20,5);(10,10) **\crv{(20,7.5)&(10,7.5)},
(20,5);(30,10) **\crv{(20,7.5)&(30,7.5)},
(10,10);(30,20) **\crv{(10,15)&(30,15)},
(30,10);(10,20) **\crv{(30,15)&(10,15)},
(50,10)*={=},
(60,0);(100,0) **@{.},
(60,20);(100,20) **@{.},
(80,0);(80,10) **@{-},
(80,10);(70,20) **\crv{(80,15)&(70,15)},
(80,10);(90,20) **\crv{(80,15)&(90,15)}
\endxy}
\eeqn\\
\beqn
\label{ess}
\vcenter{\xy
(20,35)*={s(\id_{[1]}\otimes e)};
(50,35)*={=};
(80,35)*={(e\otimes\id_{[1]})(\id_{[1]}\otimes s)(s\otimes\id_{[1]})\::};
(0,0);(40,0) **@{.},
(0,15);(40,15) **@{.},
(0,30);(40,30) **@{.},
(10,0);(25,15) **\crv{(10,7.5)&(25,7.5)},
(25,0);(10,15) **\crv{(25,7.5)&(10,7.5)},
(10,15);(10,30) **@{-},
(25,15);(25,20) **@{-},
(25,20);(20,30) **\crv{(25,25)&(20,25)},
(25,20);(30,30) **\crv{(25,25)&(30,25)},
(50,15)*={=},
(60,0);(100,0) **@{.},
(60,10);(100,10) **@{.},
(60,20);(100,20) **@{.},
(60,30);(100,30) **@{.},
(75,0);(75,5) **@{-},
(75,5);(70,10) **\crv{(75,7.5)&(70,7.5)},
(75,5);(80,10) **\crv{(75,7.5)&(80,7.5)},
(90,0);(90,10) **@{-},
(70,10);(70,20) **@{-},
(80,10);(90,20) **\crv{(80,15)&(90,15)},
(90,10);(80,20) **\crv{(90,15)&(80,15)},
(70,20);(80,30) **\crv{(70,25)&(80,25)},
(80,20);(70,30) **\crv{(80,25)&(70,25)},
(90,20);(90,30) **@{-}
\endxy}
\eeqn\\
\beqn
\label{ss}
\vcenter{\xy
(20,25)*={s^2};
(50,25)*={=};
(80,25)*={\mathrm{id}_{[2]}\::};
(0,0);(40,0) **@{.},
(0,10);(40,10) **@{.},
(0,20);(40,20) **@{.},
(10,20);(30,10) **\crv{(10,15)&(30,15)},
(30,20);(10,10) **\crv{(30,15)&(10,15)},
(10,0);(30,10) **\crv{(10,5)&(30,5)},
(30,0);(10,10) **\crv{(30,5)&(10,5)},
(50,10)*={=};
(60,0);(100,0) **@{.},
(60,20);(100,20) **@{.},
(70,0);(70,20) **@{-},
(90,0);(90,20) **@{-}
\endxy}
\eeqn\\
\beqn
\label{yangbaxter}
\vcenter{\xy
(20,35)*={(s\otimes\id_{[1]})(\id_{[1]}\otimes s)(s\otimes\id_{[1]})};
(50,35)*={=};
(80,35)*={(\id_{[1]}\otimes s)(s\otimes\id_{[1]})(\id_{[1]}\otimes s)\::};
(0,0);(40,0) **@{.},
(0,10);(40,10) **@{.},
(0,20);(40,20) **@{.},
(0,30);(40,30) **@{.},
(10,0);(20,10) **\crv{(10,5)&(20,5)},
(20,0);(10,10) **\crv{(20,5)&(10,5)},
(30,0);(30,10) **@{-},
(10,10);(10,20) **@{-},
(20,10);(30,20) **\crv{(20,15)&(30,15)},
(30,10);(20,20) **\crv{(30,15)&(20,15)},
(10,20);(20,30) **\crv{(10,25)&(20,25)},
(20,20);(10,30) **\crv{(20,25)&(10,25)},
(30,20);(30,30) **@{-},
(50,15)*={=};
(60,0);(100,0) **@{.},
(60,10);(100,10) **@{.},
(60,20);(100,20) **@{.},
(60,30);(100,30) **@{.},
(70,0);(70,10) **@{-},
(80,0);(90,10) **\crv{(80,5)&(90,5)},
(90,0);(80,10) **\crv{(90,5)&(80,5)},
(70,10);(80,20) **\crv{(70,15)&(80,15)},
(80,10);(70,20) **\crv{(80,15)&(70,15)},
(90,10);(90,20) **@{-},
(70,20);(70,30) **@{-},
(80,20);(90,30) **\crv{(80,25)&(90,25)},
(90,20);(80,30) **\crv{(90,25)&(80,25)}
\endxy}
\eeqn\\
\beqn
\label{partialc}
\vcenter{\xy
(20,25)*={c_\lambda\partial};
(50,25)*={=};
(80,25)*={\partial\otimes\partial\::};
(0,0);(40,0) **@{.},
(0,10);(40,10) **@{.},
(0,20);(40,20) **@{.},
(20,5)*+[F]{\lambda}="L";
(10,0);"L" **\crv{(10,2.5)&(20,2.5)},
(30,0);"L" **\crv{(30,2.5)&(20,2.5)},
"L";(20,10) **@{-},
(20,10);(20,15)*={\times} **@{-},
(50,10)*={=},
(60,0);(100,0) **@{.},
(60,20);(100,20) **@{.},
(70,0);(70,10)*={\times} **@{-},
(90,0);(90,10)*={\times} **@{-}
\endxy}
\eeqn\\
\beqn
\label{delred}
\vcenter{\xy
(20,25)*={c_0};
(50,25)*={=};
(80,25)*={\partial\otimes\id_{[1]}\::};
(0,0);(40,0) **@{.},
(0,20);(40,20) **@{.},
(20,10)*+[F]{0}="Z";
(20,20);"Z" **@{-},
"Z";(10,0) **\crv{(20,5)&(10,5)},
"Z";(30,0) **\crv{(20,5)&(30,5)},
(50,10)*={=},
(60,0);(100,0) **@{.},
(60,20);(100,20) **@{.},
(70,0);(70,10)*={\times} **@{-},
(90,0);(90,20) **@{-}
\endxy}
\eeqn\\
\beqn
\label{ceid}
\vcenter{\xy
(20,25)*={e\,c_\lambda};
(50,25)*={=};
(80,25)*={\id_{[1]}\::};
(0,0);(40,0) **@{.},
(0,10);(40,10) **@{.},
(0,20);(40,20) **@{.},
(20,0);(20,5) **@{-},
(20,5);(10,10) **\crv{(20,7.5)&(10,7.5)},
(20,5);(30,10) **\crv{(20,7.5)&(30,7.5)},
(20,15)*+[F]{\lambda}="L";
(10,10);"L" **\crv{(10,12.5)&(20,12.5)},
(30,10);"L" **\crv{(30,12.5)&(20,12.5)},
"L";(20,20) **@{-},
(50,10)*={=};
(60,0);(100,0) **@{.},
(60,20);(100,20) **@{.},
(80,0);(80,20) **@{-}
\endxy}
\eeqn\\
\beqn
\label{csc}
\vcenter{\xy
(20,25)*={s\,c_\lambda};
(50,25)*={=};
(80,25)*={c_{1-\lambda}\::};
(0,0);(40,0) **@{.},
(0,10);(40,10) **@{.},
(0,20);(40,20) **@{.},
(10,0);(30,10) **\crv{(10,5)&(30,5)},
(30,0);(10,10) **\crv{(30,5)&(10,5)},
(20,15)*+[F]{\lambda}="L";
(10,10);"L" **\crv{(10,12.5)&(20,12.5)},
(30,10);"L" **\crv{(30,12.5)&(20,12.5)},
"L";(20,20) **@{-},
(50,10)*={=};
(60,0);(100,0) **@{.},
(60,20);(100,20) **@{.},
(80,10)*+[F]{1-\lambda}="K";
(70,0);"K" **\crv{(70,5)&(80,5)},
(90,0);"K" **\crv{(90,5)&(80,5)},
"K";(80,20) **@{-}
\endxy}
\eeqn\\
\beqn
\label{css}
\vcenter{\xy
(20,35)*={(\id_{[1]}\otimes c_\lambda)s};
(50,35)*={=};
(80,35)*={(s\otimes \id_{[1]})(\id_{[1]}\otimes s)(c_\lambda\otimes\id_{[1]})\::};
(0,0);(40,0) **@{.},
(0,15);(40,15) **@{.},
(0,30);(40,30) **@{.},
(25,10)*+[F]{\lambda}="L";
(30,0);"L" **\crv{(30,5)&(25,5)},
(20,0);"L" **\crv{(20,5)&(25,5)},
"L";(25,15) **@{-},
(10,0);(10,15) **@{-},
(25,15);(10,30) **\crv{(25,22.5)&(10,22.5)},
(10,15);(25,30) **\crv{(10,22.5)&(25,22.5)},
(50,15)*={=},
(60,0);(100,0) **@{.},
(60,10);(100,10) **@{.},
(60,20);(100,20) **@{.},
(60,30);(100,30) **@{.},
(90,0);(90,10) **@{-},
(80,0);(70,10) **\crv{(80,5)&(70,5)},
(70,0);(80,10) **\crv{(70,5)&(80,5)},
(90,10);(80,20) **\crv{(90,15)&(80,15)},
(80,10);(90,20) **\crv{(80,15)&(90,15)},
(70,10);(70,20) **@{-},
(90,20);(90,30) **@{-},
(75,25)*+[F]{\lambda}="K";
(80,20);"K" **\crv{(80,22.5)&(75,22.5)},
(70,20);"K" **\crv{(70,22.5)&(75,22.5)},
"K";(75,30) **@{-}
\endxy}
\eeqn\\
\beqn
\label{doublec}
\vcenter{\xy
(20,35)*={(e\otimes e)(\id_{[1]}\otimes s\otimes\id_{[1]})(c_\lambda\otimes c_\lambda)};
(50,35)*={=};
(80,35)*={c_\lambda e};
(0,0);(40,0) **@{.},
(0,10);(40,10) **@{.},
(0,20);(40,20) **@{.},
(0,30);(40,30) **@{.},
(10,25)*+[F]{\lambda}="M1";
(30,25)*+[F]{\lambda}="M2";
(5,20);"M1" **\crv{(5,22.5)&(10,22.5)},
(15,20);"M1" **\crv{(15,22.5)&(10,22.5)},
"M1";(10,30) **@{-},
(25,20);"M2" **\crv{(25,22.5)&(30,22.5)},
(35,20);"M2" **\crv{(35,22.5)&(30,22.5)},
"M2";(30,30) **@{-},
(5,10);(5,20) **@{-},
(35,10);(35,20) **@{-},
(15,10);(25,20) **\crv{(15,15)&(25,15)},
(25,10);(15,20) **\crv{(25,15)&(15,15)},
(10,0);(10,5) **@{-},
(10,5);(5,10) **\crv{(10,7.5)&(5,7.5)},
(10,5);(15,10) **\crv{(10,7.5)&(15,7.5)},
(30,0);(30,5) **@{-},
(30,5);(25,10) **\crv{(30,7.5)&(25,7.5)},
(30,5);(35,10) **\crv{(30,7.5)&(35,7.5)},
(50,15)*={=},
(60,0);(100,0) **@{.},
(60,15);(100,15) **@{.},
(60,30);(100,30) **@{.},
(80,15);(80,20) **@{-},
(80,20);(70,30) **\crv{(80,25)&(70,25)},
(80,20);(90,30) **\crv{(80,25)&(90,25)},
(80,10)*+[F]{\lambda}="P";
"P";(80,15) **@{-},
(70,0);"P" **\crv{(70,5)&(80,5)},
(90,0);"P" **\crv{(90,5)&(80,5)}
\endxy}
\eeqn\\
\beqn
\label{cass}
\vcenter{\xy
(20,25)*={(c_\mu\otimes\id_{[1]})c_\lambda};
(50,25)*={=};
(80,25)*={(\id_{[1]}\otimes c_{\widetilde{\mu}})c_{\widetilde{\lambda}}\::};
(0,0);(40,0) **@{.},
(0,10);(40,10) **@{.},
(0,20);(40,20) **@{.},
(15,5)*+[F]{\mu}="M1";
(10,0);"M1" **\crv{(10,2.5)&(15,2.5)},
(20,0);"M1" **\crv{(20,2.5)&(15,2.5)},
"M1";(15,10) **@{-},
(30,0);(30,10) **@{-},
(22.5,15)*+[F]{\lambda}="L1";
(15,10);"L1" **\crv{(15,12.5)&(22.5,12.5)},
(30,10);"L1" **\crv{(30,12.5)&(22.5,12.5)},
"L1";(22.5,20) **@{-},
(50,10)*={=},
(60,0);(100,0) **@{.},
(60,10);(100,10) **@{.},
(60,20);(100,20) **@{.},
(70,0);(70,10) **@{-},
(85,5)*+[F]{\widetilde{\mu}}="M2";
(80,0);"M2" **\crv{(80,2.5)&(85,2.5)},
(90,0);"M2" **\crv{(90,2.5)&(85,2.5)},
"M2";(85,10) **@{-},
(77.5,15)*+[F]{\widetilde{\lambda}}="K";
(70,10);"K" **\crv{(70,12.5)&(77.5,12.5)},
(85,10);"K" **\crv{(85,12.5)&(77.5,12.5)},
"K";(77.5,20) **@{-},
\endxy}
\eeqn\\
using the abbreviations
\beq
\widetilde{\lambda}=\lambda\mu,\qquad\widetilde{\mu}=\left\{\begin{array}{cl}\lambda\frac{1-\mu}{1-\lambda\mu}&\textrm{ if }\lambda\mu\neq 1\\\textrm{arbitrary}&\textrm{ if }\lambda=\mu=1.\end{array}\right.
\eeq
\end{defn}

Hence, a morphism in $\mathtt{FinStoMap'}$ is represented by a vertical juxtaposition of horizontal juxtapositions of generators and identity morphisms such that the strands match. Two such diagrams describe the same morphism if and only if there is a sequence of steps of the form\re{eass}--(\ref{cass}) transforming the two diagrams into each other. The way to think of a diagrammatic representation of a morphism in $\mathtt{FinStoMap'}$ is as a probabilistic map $[m]\ra[n]$, where the image of $j\in[m]$ can be obtained by following the $j$th input strand downwards, such that at an occurence of some $c_\lambda$ one branches to the left with probability $\lambda$ and branches to the right with probability $1-\lambda$. One can check easily that the defining relations of $\mathtt{FinStoMap'}$ are consistent with this interpretation.\\

\textsc{Remark.}
\begin{compactenum}
\item By combining\re{delred} with\re{ceid} and\re{csc}, we obtain two additional useful equations:\\\\
\beqn
\label{edelisid}
\vcenter{\xy
(20,25)*={e(\partial\otimes\mathrm{id}_{[1]})};
(50,25)*={=};
(80,25)*={\mathrm{id}_{[1]}\::};
(0,0);(40,0) **@{.},
(0,10);(40,10) **@{.},
(0,20);(40,20) **@{.},
(30,20);(30,10) **@{-},
(10,15)*={\times};(10,10) **@{-},
(20,5);(20,0) **@{-},
(20,5);(10,10) **\crv{(20,7.5)&(10,7.5)},
(20,5);(30,10) **\crv{(20,7.5)&(30,7.5)},
(50,10)*={=};
(60,0);(100,0) **@{.},
(60,20);(100,20) **@{.},
(80,0);(80,20) **@{-},
\endxy}
\eeqn\\
\beqn
\label{psp}
\vcenter{\xy
(20,25)*={s(\partial\otimes\id_{[1]})};
(50,25)*={=};
(80,25)*={\id_{[1]}\otimes\partial\::};
(0,0);(40,0) **@{.},
(0,10);(40,10) **@{.},
(0,20);(40,20) **@{.},
(10,0);(30,10) **\crv{(10,5)&(30,5)},
(30,0);(10,10) **\crv{(30,5)&(10,5)},
(30,10);(30,20) **@{-},
(10,10);(10,15)*={\times} **@{-},
(50,10)*={=};
(60,0);(100,0) **@{.},
(60,20);(100,20) **@{.},
(70,0);(70,20) **@{-},
(90,0);(90,10)*={\times} **@{-}
\endxy}
\eeqn\\
As proven in~\cite{Mas}, none of the equations\re{eass}--(\ref{yangbaxter}),\re{edelisid},\re{psp} which form an analogous presentation of $\mathtt{FinMap}$ (where the generators $c_\lambda$ are not present) is implied by the other six.
\item As already noted in~\cite{Laf}, the equations\re{ess},\re{edelisid} and\re{psp} imply their mirror images by use of\re{ss} and\re{ese}. The same holds true for\re{css}.
\item As can be seen from the relation\re{delred}, the generator $\partial$ is redundant for all morphisms $f:[m]\ra[n]$ with $m\geq 1$. Hence its only function is to turn $[0]$ into an initial object in $\mathtt{FinStoMap}$, as without $\partial$ there could be no morphism from $[0]$ to any other object.\\
\end{compactenum}

Taking the strict monoidal functor $F:\mathtt{FinStoMap'}\ra\mathtt{FinStoMap}$ to be the identity on objects, the assignments
\beqa
F(\partial)&\equiv&\left(\begin{array}{c}\end{array}\right)  \::\:  [0]\lra [1]\\\\
F(e)&\equiv& \left(\begin{array}{cc}1 & 1\end{array}\right) \::\:  [2]\lra [1]\\\\
F(s)&\equiv&\left(\begin{array}{cc}0&1\\1&0\end{array}\right) \::\:  [2]\lra [2]\\\\
F(c_\lambda)&\equiv&\left(\begin{array}{c}\lambda\\1-\lambda\end{array}\right) \::\: [1]\lra [2]
\eeqa
preserve the relations and hence uniquely define $F$. The motivation for these definitions is that they exactly match the interpretations of the generators of $\mathtt{FinStoMap'}$ as the corresponding probabilistic maps. When a stochastic matrix $A$ has a preimage $F^{-1}(A)$ in $\mathtt{FinStoMap'}$, this preimage then provides a possible way to turn the blank rectangle of the ``black box''\re{blackbox} into a concrete representation of strands branching, crossing, coalescing, and newly emerging.

The series of intermediate results following now will culminate in theorem~\ref{FinStoMapPres} stating that the functor $F$ is in fact an isomorphism of strict monoidal categories.

\begin{lem}
For $n\geq 1$, every morphism $f\in\mathtt{FinStoMap'}([1],[n])$ can be written in the form
\beqn
\label{fasc}
f=(\id_{[n-2]}\otimes c_{\lambda_{n-1}})\cdots(\id_{[1]}\otimes c_{\lambda_2})c_{\lambda_1}
\eeqn
with numbers $\lambda_j\in[0,1]$. The image $F(f)$ is a stochastic matrix
\beq
F(f)=\left(\begin{array}{c}\mu_1\\\vdots\\\mu_{n-1}\\\eta_n\end{array}\right)
\eeq
with entries
\beqn
\label{smsystem}
\mu_j=\lambda_j(1-\lambda_{j-1})\cdots(1-\lambda_1),\quad j=1,\ldots n-1;\qquad\eta_n=(1-\lambda_{n-1})\cdots(1-\lambda_1).
\eeqn
\end{lem}

It is understood that\re{fasc} degerenates to the empty product when $n=1$, i.e. the statement is that $f=\id_{[1]}$ in this case.

\begin{proof}
First, it will be shown that any such $f$ can be written without using the generators $\partial$, $e$, or $s$. For $\partial$, this is clear by the relation\re{delred}. Then we may write $f$ as a product of terms of the form $\id_{[\cdot]}\otimes e\otimes\id_{[\cdot]}$, $\id_{[\cdot]}\otimes s\otimes\id_{[\cdot]}$, and $\id_{[\cdot]}\otimes c_\lambda\otimes\id_{[\cdot]}$. Now consider the rightmost term in this product which contains a generator $e$ or $s$ and hence has the form $\id_{[k]}\otimes e\otimes\id_{[l]}$ or $\id_{[k]}\otimes s\otimes\id_{[l]}$. Such a factor has $k+l+2$ input strands. Since $f$ itself only has a single input strand, there have to be exactly $k+l+1$ factors to the right of it, each being of the form $\id_{[\cdot]}\otimes c_\lambda\otimes\id_{[\cdot]}$. Hence by repeated application of deformed parametric associativity\re{cass}, we can write $f$ in such a form that the factor immediately succeeding the $\id_{[k]}\otimes e\otimes\id_{[l]}$ or $\id_{[k]}\otimes s\otimes\id_{[l]}$ has the form $\id_{[k]}\otimes c_\lambda\otimes\id_{[l]}$. Then an application of the relation\re{ceid} or\re{csc} removes the occurence of the unwanted generator $e$ or $s$. This procedure now can be applied repeatedly until all occurences of $e$ and $s$ are removed. We now have a representation of $f$ with exactly $n-1$ factors of the form $\id_{[\cdot]}\otimes c_\lambda\otimes\id_{[\cdot]}$ and containing no other generators.

Second, again by repeated application of deformed parametric associativity\re{cass}, $f$ then can be brought into the form whose existence was asserted.

For the second assertion, apply induction on $n$. For $n=1$, there is nothing to prove. Taking the assertion for $n$ as the induction assumption, we get for the case of $n+1$ that
\beqa
&F\left((\id_{[n-1]}\otimes c_{\lambda_n})\cdots(\id_{[1]}\otimes c_{\lambda_2})c_{\lambda_1}\right)\\\\
&=\left(\mathbbm{1}_{n-1}\otimes F(c_\lambda)\right)F\left((\id_{[n-2]}\otimes c_{\lambda_{n-1}})\cdots(\id_{[1]}\otimes c_{\lambda_2})c_{\lambda_1}\right)\\
&=\left(\begin{array}{cc}\mathbbm{1}_{n-1}&0\\0&\lambda_n\\0&1-\lambda_n\end{array}\right)\left(\begin{array}{c}\mu_1\\\vdots\\\mu_{n-1}\\\eta_n\end{array}\right)=\left(\begin{array}{c}\mu_1\\\vdots\\\mu_{n-1}\\\lambda_n\eta_n\\(1-\lambda_n)\eta_n\end{array}\right)=\left(\begin{array}{c}\mu_1\\\vdots\\\mu_n\\\eta_{n+1}\end{array}\right).
\eeqa
\end{proof}

Now we can use this result to prove that $F$ is bijective on those morphism sets that have the object $[1]$ as their domain. The rest of this article then will be devoted to proving that a morphism in $\mathtt{FinStoMap'}([m],[n])$ can be decomposed into $m$ morphisms in $\mathtt{FinStoMap'}([1],[n])$ in a way that is compatible with decomposing a stochastic matrix in $\mathtt{FinStoMap}([m],[n])$ into its $m$ columns in $\mathtt{FinStoMap}([1],[n])$.

\begin{prop}
\label{Fsingleiso}
For every $n\in\N_0$, the map $F([1],[n]):\mathtt{FinStoMap'}([1],[n])\ra\mathtt{FinStoMap}([1],[n])$ is bijective.
\end{prop}

\begin{proof}
This is clear for $n=0$, as both $\mathtt{FinStoMap'}([1],[0])$ and $\mathtt{FinStoMap}([1],[0])$ are empty. For $n\geq 1$, suppose that we have a single-column stochastic matrix
\beq
A_n=\left(\begin{array}{c}\mu_1\\\vdots\\\mu_{n-1}\\\eta_n\end{array}\right)
\eeq
with entries $\mu_j\geq 0$, $\eta_n\geq 0$ satisfying $\eta_n=1-\sum_j\mu_j$. This matrix has a preimage under $F$ of the form\re{fasc} if we can solve the system\re{smsystem} for appropriate $\lambda_j\in[0,1]$. An explicit solution is given by
\beq
\label{solution}
\lambda_j=\frac{\mu_j}{1-\sum_{k=1}^{j-1}\mu_k},\quad j=1,\ldots,n-1
\eeq
with the convention that $0/0$ may be an arbitrary value in $[0,1]$. It can be verified by direct calculation that this solves\re{smsystem}. As for uniqueness, note that the system of equations\re{smsystem} can also be solved for the $\lambda_j$ recursively starting with $\lambda_1=\mu_1$, as long as we never have $\lambda_j=1$ for some $j$. In this exceptional case, we can take $\lambda_k$ to be arbitrary for $k>j$. Hence the proof is complete if we can show that we get the same morphism in $\mathtt{FinStoMap'}$ no matter which choice of $\lambda_k$, $k>j$, we make in this case. This follows from repeated application of the equation
\beq
\vcenter{\xy
(20,25)*={(\id_{[1]}\otimes c_\lambda)c_1};
(50,25)*={=};
(80,25)*={(\id_{[1]}\otimes c_1)c_1};
(0,0);(40,0) **@{.},
(0,10);(40,10) **@{.},
(0,20);(40,20) **@{.},
(25,5)*+[F]{\lambda}="M1";
(20,0);"M1" **\crv{(20,2.5)&(25,2.5)},
(30,0);"M1" **\crv{(30,2.5)&(25,2.5)},
"M1";(25,10) **@{-},
(10,0);(10,10) **@{-},
(17.5,15)*+[F]{1}="L1";
(25,10);"L1" **\crv{(25,12.5)&(17.5,12.5)},
(10,10);"L1" **\crv{(10,12.5)&(17.5,12.5)},
"L1";(17.5,20) **@{-},
(50,10)*={=},
(60,0);(100,0) **@{.},
(60,10);(100,10) **@{.},
(60,20);(100,20) **@{.},
(70,0);(70,10) **@{-},
(85,5)*+[F]{1}="M2";
(80,0);"M2" **\crv{(80,2.5)&(85,2.5)},
(90,0);"M2" **\crv{(90,2.5)&(85,2.5)},
"M2";(85,10) **@{-},
(77.5,15)*+[F]{1}="K";
(70,10);"K" **\crv{(70,12.5)&(77.5,12.5)},
(85,10);"K" **\crv{(85,12.5)&(77.5,12.5)},
"K";(77.5,20) **@{-},
\endxy}
\eeq
which is a consequence of deformed parametric associativity\re{cass}.
\end{proof}

Now what we have to do is to set up a bijection between $\mathtt{FinStoMap'}([m],[n])$ and $\mathtt{FinStoMap'}([1],[n])^m$, such that this bijection corresponds under $F$ to decomposing a stochastic matrix into its columns. The hardest part of this is to specify how to obtain a morphism in $\mathtt{FinStoMap'}([m],[n])$, given an $m$-tuple of morphisms in $\mathtt{FinStoMap'}([1],[n])$. Taking the tensor product of the elements of the original $m$-tuple produces a morphism $[m]\ra[mn]$. Then by the yet to be defined family of coalescing maps $p_n^m:[mn]\ra[n]$, we obtain the composition $[m]\ra[mn]\ra[n]$, which is the desired element in $\mathtt{FinStoMap'}([m],[n])$.

Before this family of coalescing maps $p_n^m$ can be introduced, it is necessary to study another family of particular morphisms in $\mathtt{FinStoMap'}$ and to prove some formulas about them. The ``cyclic permutation'' morphisms $z_n:[n]\ra [n]$ are defined recursively via
\beqn
\label{defz}
z_1\equiv\id_{[1]};\qquad z_{n+1}\equiv(\id_{[n-1]}\otimes s)(z_n\otimes\id_{[1]}),\quad n\geq 1.
\eeqn
The morphism $z_n$ can be thought of as a permutation of the $n$ strands which turns the leftmost strand into the rightmost strand while keeping the order of the other strands fixed. As we will see now, this interpretation is confirmed by the image of $z_n$ in $\mathtt{FinStoMap}$.

\begin{lem}
The functor $F$ maps $z_n$ to the permutation matrix which turns the leftmost strand into the rightmost strand while keeping the order of the other strands fixed:
\beqn
\label{zmatrix}
F(z_n)=\left(\begin{array}{cc}0&\mathbbm{1}_{n-1}\\1&0\end{array}\right)
\eeqn
\end{lem}

\begin{proof}
Again induction on $n$. The case $n=1$ is clear. Then,
\beqa
&F(z_{n+1})=F(\id_{[n-1]}\otimes s)F(z_n\otimes\id_{[1]})=\left(\begin{array}{ccc}\mathbbm{1}_{n-1}&0&0\\0&0&1\\0&1&0\end{array}\right)\left(\begin{array}{ccc}0&\mathbbm{1}_{n-1}&0\\1&0&0\\0&0&1\end{array}\right)\\
&=\left(\begin{array}{ccc}0&\mathbbm{1}_{n-1}&0\\0&0&1\\1&0&0\end{array}\right)=\left(\begin{array}{cc}0&\mathbbm{1}_n\\1&0\end{array}\right)
\eeqa
\end{proof}

The results of the next lemma are immediate if one knows that the generator $s$ and the relations\re{ss},\re{yangbaxter} form a presentation of the strict monoidal category of invertible maps on finite cardinals. For the sake of completeness, we give an independent proof here.

\begin{lem}
The cyclic permutation morphisms $z_n$ are invertible and satisfy the following equations:
\begin{compactenum}
\item For any integer $n\geq 1$,
\beqn
\label{otherdefz}
z_{n+1}=(\id_{[1]}\otimes z_{n})(s\otimes\id_{[n-1]}).
\eeqn
\item For any integer $n\geq 1$,
\beq
z_n\otimes z_n=(\id_{[n-1]}\otimes z_{n+1})(z_{n+1}\otimes\id_{[n-1]}).
\eeq
\item For any integer $n\geq 1$,
\beqn
\label{zeq}
(z_n^{-1}\otimes\id_{[n]})(\id_{[n-1]}\otimes z_{n+1})=(\id_{[n]}\otimes z_n)(z_{n+1}^{-1}\otimes\id_{[n-1]}).
\eeqn
\item For any integer $n\geq 0$,
\beqn
\label{zpartial}
z_{n+1}(\partial\otimes\id_{[n]})=\id_{[n]}\otimes\partial.
\eeqn
\end{compactenum}
\end{lem}

\begin{proof}
Invertibility is clear as $z_n$ is defined as a composition of invertible morphisms. All the following proofs use induction on $n$. 
\begin{compactenum}
\item
Trivial for $n=1$, while the induction step is
\beqa
&z_{n+2}\stackrel{\re{defz}}{=}(\id_{[n]}\otimes s)(z_{n+1}\otimes\id_{[1]})\stackrel{\textrm{i.a.}}{=}(\id_{[n]}\otimes s)(\id_{[1]}\otimes z_n\otimes\id_{[1]})(s\otimes\id_{[n]})\\
&\stackrel{\re{defz}}{=}(\id_{[1]}\otimes z_{n+1})(s\otimes\id_{[n]})
\eeqa
\item The case $n=1$ states $\id_{[1]}\otimes\id_{[1]}=ss$, which is\re{ss}. The following calculation proves the assertion for $n+1$ assuming its validity for $n$:
\beqa
&z_{n+1}\otimes z_{n+1}=(\id_{[n+1]}\otimes z_{n+1})(z_{n+1}\otimes\id_{[n+1]})\\
&\stackrel{\re{defz},\re{otherdefz}}{=}(\id_{[2n]}\otimes s)(\id_{[n+1]}\otimes z_n\otimes\id_{[1]})(\id_{[1]}\otimes z_n\otimes\id_{[n+1]})(s\otimes\id_{[2n]})\\
&\stackrel{\textrm{i.a.}}{=}(\id_{[2n]}\otimes s)(\id_{[n]}\otimes z_{n+1}\otimes\id_{[1]})(\id_{[1]}\otimes z_{n+1}\otimes\id_{[n]})(s\otimes\id_{[2n]})\\
&\stackrel{\re{defz},\re{otherdefz}}{=}(\id_{[n]}\otimes z_{n+2})(z_{n+2}\otimes\id_{[n]})
\eeqa
\item This is the previous equation in a different form.
\item The statement is vacuous for $n=0$. The induction step is
\beqa
&z_{n+2}(\partial\otimes\id_{[n+1]})\stackrel{\re{defz}}{=}(\id_{[n]}\otimes s)(z_{n+1}\otimes\id_{[1]})(\partial\otimes\id_{[n+1]})\\
&\stackrel{\textrm{i.a.}}{=}(\id_{[n]}\otimes s)(\id_{[n]}\otimes\partial\otimes\id_{[1]})\stackrel{\re{psp}}{=}\id_{[n+1]}\otimes\partial
\eeqa
\end{compactenum}
\end{proof}

The next lemma then uses\re{otherdefz} and some of the relations in $\mathtt{FinStoMap'}$ to study how the $z_n$ behave with respect to arbitrary morphisms in $\mathtt{FinStoMap'}$.

\begin{lem}
\label{zcomm}
For $f\in\mathtt{FinStoMap'}([m],[n])$, we have
\beqn
\label{zcommeq}
z_{n+1}(\id_{[1]}\otimes f)=(f\otimes\id_{[1]})z_{m+1}.
\eeqn
\end{lem}

\begin{proof}
This will be done in the following three steps:
\begin{compactenum}
\item It holds for $f=\partial$, $e$, $s$ and all $c_\lambda$.
\item If it holds for $f$, then it also holds for any $\id_{[k]}\otimes f\otimes\id_{[l]}$.
\item If it holds for $f_1:[m]\ra[n]$ and $f_2:[n]\ra[q]$, then it also holds for $f_2f_1:[m]\ra[q]$.
\end{compactenum}
This then covers all cases as every morphism is a composition of tensor products of generators and identity morphisms.
\begin{compactenum}
\item For $f=\partial$, this is\re{psp}. For $f=e$, it is\re{ess}. For $f=s$ itself, this is the Yang-Baxter relation\re{yangbaxter}, while for $c_\lambda$ it is\re{css}. 
\item It is sufficient to prove this for the cases $k=0$, $l=1$ and $k=1$, $l=0$, as all other cases then follow by induction. For the first of these, this is the calculation
\beqa
&z_{n+2}(\id_{[1]}\otimes f\otimes\id_{[1]})=(\id_{[n]}\otimes s)(z_{n+1}\otimes\id_{[1]})(\id_{[1]}\otimes f\otimes\id_{[1]})\\
&=(\id_{[n]}\otimes s)(f\otimes\id_{[2]})(z_{m+1}\otimes\id_{[1]})=(f\otimes\id_{[2]})(\id_{[m]}\otimes s)(z_{m+1}\otimes\id_{[1]})\\
&=(f\otimes\id_{[2]})z_{m+2}
\eeqa
while the second case works similarly using\re{otherdefz}.
\item Direct calculation:
\beqa
&z_{q+1}(\id_{[1]}\otimes f_2f_1)=z_{q+1}(\id_{[1]}\otimes f_2)(\id_{[1]}\otimes f_1)=(f_2\otimes\id_{[1]})z_{n+1}(\id_{[1]}\otimes f_1)\\
&=(f_2\otimes\id_{[1]})(f_1\otimes\id_{[1]})z_{m+1}=(f_2f_1\otimes\id_{[1]})z_{m+1}
\eeqa
\end{compactenum}
\end{proof}

Now the coalescing morphisms $p_n^m:[mn]\ra[n]$ can be introduced. $p_n^m$ coalesces $m$ copies of a group of $n$ strands into a single group of $n$ strands and can be defined recursively by (with $n\geq 0$, $m\geq 2$)
\beqn
\label{defp}
p^2_0\equiv\id_{[0]},\qquad p^2_{n+1}\equiv (p^2_n\otimes e)(\id_{[n]}\otimes z_{n+2}),\qquad p_n^{m+1}\equiv p^2_n(p^m_n\otimes\id_{[n]})
\eeqn
The interpretation of $p_m^n$ as coalescing strands is confirmed by its image in $\mathtt{FinStoMap}$:

\begin{lem}
\label{pcopies}
For integers $m\geq 2$ and $n\geq 0$,
\beq
F(p^m_n)=\underset{m\textrm{ copies}}{\underbrace{\left(\mathbbm{1}_n\cdots\mathbbm{1}_n\right)}}.
\eeq
\end{lem}

\begin{proof}
First, induction on $n$ for $m=2$:
\beqa
&F(p_{n+1}^2)=F(p_n^2\otimes e)F(\id_{[n]}\otimes z_{n+2})=\left(\begin{array}{cccc}\mathbbm{1}_n&\mathbbm{1}_n&0&0\\0&0&1&1\end{array}\right)\left(\begin{array}{cccc}\mathbbm{1}_n&0&0&0\\0&0&\mathbbm{1}_n&0\\0&0&0&1\\0&1&0&0\end{array}\right)\\
&=\left(\begin{array}{cccc}\mathbbm{1}_n&0&\mathbbm{1}_n&0\\0&1&0&1\end{array}\right)=\left(\begin{array}{cc}\mathbbm{1}_{n+1}&\mathbbm{1}_{n+1}\end{array}\right)
\eeqa
Then, induction on $m$ for fixed $n$:
\beqa
&F(p^{m+1}_n)=F(p_n^2)F(p_n^m\otimes\id_{[n]})\\\\
&=\left(\begin{array}{cc}\mathbbm{1}_n&\mathbbm{1}_n\end{array}\right)\left(\begin{array}{cccc}\mathbbm{1}_n&\cdots&\mathbbm{1}_n&0\\0&\cdots&0&\mathbbm{1}_n\end{array}\right)=\left(\begin{array}{cccc}\mathbbm{1}_n&\cdots&\mathbbm{1}_n&\mathbbm{1}_n\end{array}\right)
\eeqa
\end{proof}

Similar to\re{otherdefz} for the $z_n$'s, it will be necessary to also have another expression for $p_{n+1}^2$ in terms of $p_n^2$.

\begin{lem}
For integer $n\geq 0$,
\beqn
\label{otherdefp}
p_{n+1}^2=(e\otimes p_n^2)(z_{n+2}^{-1}\otimes\id_{[n]}).
\eeqn
\end{lem}

\begin{proof}
Induction on $n$. The statement is trivial for $n=0$. The induction step is
\beqa
&p_{n+2}^2\stackrel{\re{defp}}{=}(p_{n+1}^2\otimes e)(\id_{[n+1]}\otimes z_{n+3})\stackrel{\textrm{i.a.}}{=}(e\otimes p_n^2\otimes e)(z_{n+2}^{-1}\otimes\id_{[n+2]})(\id_{[n+1]}\otimes z_{n+3})\\
&\stackrel{\re{zeq}}{=}(e\otimes p_n^2\otimes e)(\id_{[n+2]}\otimes z_{n+2})(z_{n+3}^{-1}\otimes\id_{[n+1]})\stackrel{\re{defp}}{=}(e\otimes p_{n+1}^2)(z_{n+3}^{-1}\otimes\id_{[n]})
\eeqa
\end{proof}

The next lemma is the most important one. Similar to what lemma~\ref{zcomm} did for the $z_n$'s, it shows that the $p_n^m$'s commute with arbitrary morphisms in $\mathtt{FinStoMap'}$ in a certain way.

\begin{lem}
\label{pcomm}
For any $f:[m]\ra[n]$ and any integer $k\geq 2$, we have
\beq
fp^k_m=p^k_nf^{\otimes k}
\eeq
\end{lem}

\begin{proof}
Consider the case $k=2$ first. This then uses exactly the same three steps as the proof of lemma~\ref{zcomm} did.
\begin{compactenum}
\item We have $p_1^2=es=e$, and hence 
\beq
p_2^2=(e\otimes e)(\id_{[2]}\otimes s)(\id_{[1]}\otimes s\otimes\id_{[1]})=(e\otimes e)(\id_{[1]}\otimes s\otimes \id_{[1]}).
\eeq
For $f=\partial$, the assertion $\partial=e(\partial\otimes\partial)$ then directly follows from\re{edelisid}. For $f=e$, we need\re{ese} together with several applications of\re{eass}. For $f=s$, the calculation uses\re{ss} as well as several applications of\re{ess} and its mirror image. Finally, for $f=c_\lambda$, this is\re{doublec}.
\item Straightforward calculation employing lemma~\ref{zcomm}:
\beqa
&p_{n+1}^2(f\otimes\id_{[1]}\otimes f\otimes\id_{[1]})=(p_n^2\otimes e)(\id_{[n]}\otimes z_{n+2})(f\otimes\id_{[1]}\otimes f\otimes\id_{[1]})\\
&\stackrel{\re{zcommeq}}{=}(p_n^2\otimes e)(f\otimes f\otimes\id_{[2]})(\id_{[m]}\otimes z_{m+2})\stackrel{\textrm{i.a.}}{=}(f\otimes\id_{[1]})(p_m^2\otimes e)(\id_{[m]}\otimes z_{m+2})\\
&=(f\otimes\id_{[1]})p_{m+1}^2
\eeqa
as well as
\beqa
&p_{n+1}^2(\id_{[1]}\otimes f\otimes\id_{[1]}\otimes f)\stackrel{\re{otherdefp}}{=}(e\otimes p_n^2)(z_{n+2}^{-1}\otimes\id_{[n]})(\id_{[1]}\otimes f\otimes\id_{[1]}\otimes f)\\
&\stackrel{\re{zcommeq}}{=}(e\otimes p_n^2)(\id_{[2]}\otimes f\otimes f)(z_{m+2}^{-1}\otimes\id_{[m]})\stackrel{\textrm{i.a.}}{=}(\id_{[1]}\otimes f)(e\otimes p_m^2)(z_{m+2}^{-1}\otimes\id_{[m]})\\
&\stackrel{\re{otherdefp}}{=}(\id_{[1]}\otimes f)p_{m+1}^2.
\eeqa
\item Again the same simple calculation as in the proof of lemma~\ref{zcomm} (also using the same notation):
\beq
f_2f_1p^2_m=f_2p^2_nf_1^{\otimes 2}=p^2_qf_2^{\otimes 2}f_1^{\otimes 2}=p^2_q(f_2f_1)^{\otimes 2}
\eeq
\end{compactenum}
For general $k$, the statement is an easy consequence of the $k=2$ case and the definition\re{defp}. Upon induction on $k$,
\beqa
&fp^{k+1}_m=fp^2_m(p^k_m\otimes\id_{[m]})=p^2_n(fp^k_m\otimes f)=p^2_n(p^k_nf^{\otimes k}\otimes f)\\
&=p^2_n(p^k_n\otimes\id_{[n]})f^{\otimes(k+1)}=p^{k+1}_nf^{\otimes(k+1)}.
\eeqa
\end{proof}

\begin{lem}
For all integers $n\geq m\geq 0$,
\beqn
\label{p2del}
p_n^2(\id_{[m]}\otimes\partial^{\otimes n}\otimes\id_{[n-m]})=\id_{[n]}.
\eeqn
\end{lem}

\begin{proof}
Induction on $n$. For $n=0$, there is nothing to prove, hence proceed to the induction step and let us show that the equation holds for $n+1$ if it holds for $n$. Consider the case $m
\leq n$ first. Then the assertion follows as in
\beqa
&p_{n+1}^2(\id_{[m]}\otimes\partial^{\otimes(n+1)}\otimes\id_{[n+1-m]})=(p_n^2\otimes e)(\id_{[n]}\otimes z_{n+2})(\id_{[m]}\otimes\partial^{\otimes(n+1)}\otimes\id_{[n+1-m]})\\
&=(p_n^2\otimes e)\left[\id_{[m]}\otimes\partial^{\otimes(n-m)}\otimes z_{n+2}(\partial\otimes\id_{[n+1]})(\partial^{\otimes m}\otimes\id_{[n+1-m]})\right]\\
&\stackrel{\re{zpartial}}{=}(p_n^2\otimes e)\left[\id_{[m]}\otimes\partial^{\otimes(n-m)}\otimes(\id_{[n+1]}\otimes\partial)(\partial^{\otimes m}\otimes\id_{[n+1-m]})\right]\\
&=(p_n^2\otimes e)(\id_{[m]}\otimes\partial^{\otimes n}\otimes\id_{[n+1-m]}\otimes\partial)\underset{\re{psp}}{\overset{\textrm{i.a.}}{=}}\id_{[n]}\otimes\id_{[1]}.
\eeqa
In the case that $m=n+1$, we can use\re{otherdefp} to complete the induction step:
\beqa
&p_{n+1}^2(\id_{[n+1]}\otimes\partial^{\otimes(n+1)})=(e\otimes p_n^2)(z_{n+2}^{-1}\otimes\id_{[n]})(\id_{[n+1]}\otimes\partial^{\otimes(n+1)})\\
&\stackrel{\re{zpartial}}{=}(e\otimes p_n^2)(\partial\otimes\id_{[n+1]}\otimes\partial^{\otimes n})\underset{\re{psp}}{\overset{\textrm{i.a.}}{=}}\id_{[1]}\otimes\id_{[n]}
\eeqa
\end{proof}

Finally, a last class of morphisms in $\mathtt{FinStoMap'}$ needs to be introduced. The single-strand inclusion
\beq
\iota_j^n\equiv \partial_{j-1}\otimes\id_{[1]}\otimes\partial_{n-j}
\eeq
is a morphism $[1]\ra[n]$ which maps a single input strand to the $j$th of $n$
output strands. The composition $f\iota_j^m$ for some morphism $f:[m]\ra[n]$
then is the morphism $[1]\ra[n]$ which should be interpreted as the $j$th
``column'' of $f$. We now have to prove that these inclusion morphism are compatible with the coalescing morphisms in the expected way:

\begin{lem}
\label{piotas}
For all integers $n\geq m\geq 2$,
\beq
p_n^m(\iota_1^n\otimes\cdots\otimes\iota_m^n)=\id_{[m]}\otimes\partial^{\otimes(n-m)}.
\eeq
\end{lem}

\begin{proof}
For $m=2$, apply induction on $n$. The case $n=2$ is a direct calculation using $p_2^2=(e\otimes e)(\id_{[1]}\otimes s\otimes\id_{[1]})$ together with the equations\re{edelisid} and\re{psp}. The induction step is
\beqa
&p_{n+1}^2(\iota_1^{n+1}\otimes\iota_2^{n+1})=(p_n^2\otimes e)(\id_{[n]}\otimes z_{n+2})(\iota_1^n\otimes\partial\otimes\iota_2^n\otimes\partial)\\
&\stackrel{\re{zpartial}}{=}(p_n^2\otimes e)(\iota_1^n\otimes\iota_2^n\otimes\partial\otimes\partial)\underset{\re{edelisid}}{\overset{\textrm{i.a.}}{=}}\id_{[2]}\otimes\partial^{\otimes(n-2)}\otimes\partial=\id_{[2]}\otimes\partial^{\otimes(n-1)}.
\eeqa
Finally, we use induction on $m$:
\beqa
&p_n^{m+1}(\iota_1^n\otimes\cdots\otimes\iota_m^n\otimes\iota_{m+1}^n)=p_n^2(p_n^m\otimes\id_{[n]})(\iota_1^n\otimes\cdots\otimes\iota_m^n\otimes\iota_{m+1}^n)\\
&\stackrel{\textrm{i.a.}}{=}p_n^2(\id_{[m]}\otimes\partial^{\otimes(n-m)}\otimes\iota_{m+1}^n)=p_n^2(\id_{[m]}\otimes\partial^{\otimes n}\otimes\id_{[1]}\otimes\partial^{\otimes(n-m-1)})\\
&=p_n^2(\id_{[m]}\otimes\partial^{\otimes n}\otimes\id_{[n-m]})(\id_{[m+1]}\otimes\partial^{\otimes(n-m-1)})\stackrel{\re{p2del}}{=}\id_{[m+1]}\otimes\partial^{\otimes(n-m-1)}
\eeqa
\end{proof}

In order for the following two propositions to make sense also in the cases $m=0$ and $m=1$, let us set $p^1_n=\id_{[n]}$ and $p^0_n=\partial^{\otimes n}$. Then lemma~\ref{pcopies} immediately extends to these cases. Now after the preparations are done, we can prove the desired decomposition of a morphism in $\mathtt{FinStoMap}$ into its ``columns''.

\begin{prop}
\label{finalprop}
For any morphism $f:[m]\ra[n]$ in $\mathtt{FinStoMap'}$,
\beq
f=p_n^m(f\iota_1^m\otimes\ldots\otimes f\iota_m^m)
\eeq
\end{prop}

\begin{proof}
For $m=1$, the statement is trivial. For $m\geq 2$, this is an immediate consequence of the two lemmas~\ref{pcomm} and~\ref{piotas}. It remains to consider the degenerate case $m=0$, where the equation asserts that $f=\partial^{\otimes n}$. But this in turn follows from repeated applications of\re{partialc},\re{edelisid} and\re{psp}.
\end{proof}

That this decomposition indeed corresponds to the decomposition of a stochastic matrix into its columns is then expressed by the next proposition. 

\begin{prop}
\label{matrixdecomp}
For any stochastic matrix $A:[m]\ra[n]$, we have
\beq
A=F(p_n^m)\left(AF(\iota_1^m)\otimes\ldots\otimes AF(\iota_m^m)\right)
\eeq
\end{prop}

\begin{proof}
By definition, $F(\iota_j^m)$ is the single-column matrix with a $1$ as the $j$th entry and zeros otherwise. Hence, $A_j\equiv AF(\iota_j^m)$ is simply the $j$th column of $A$. Consequently,
\beqa
&F(p_n^m)\left(AF(\iota_1^m)\otimes\ldots\otimes AF(\iota_m^m)\right)=\left(\begin{array}{ccc}\mathbbm{1}_m&\cdots&\mathbbm{1}_m\end{array}\right)\left(\begin{array}{ccc}A_1&&0\\&\ddots\\0&&A_m\end{array}\right)\\
&=\left(\begin{array}{ccc}A_1&\cdots&A_m\end{array}\right)=A
\eeqa
\end{proof}

\begin{thm}
\label{FinStoMapPres}
The functor $F:\mathtt{FinStoMap'}\ra\mathtt{FinStoMap}$ is an isomorphism of strict monoidal categories.
\end{thm}

\begin{proof}
The two previous propositions show that a morphism $f\in\mathtt{FinStoMap'}([m],[n])$ or $A\in\mathtt{FinStoMap}([m],[n])$ is uniquely determined by an $m$-tuple of morphisms $(f\iota_j)_j$ in $\mathtt{FinStoMap'}([1],[n])$ or $\left(AF(\iota_j)\right)_j$ in $\mathtt{FinStoMap}([1],[n])$, respectively. This is expressed by the two horizontal bijections in the diagram
\beq
\xymatrix{{\mathtt{FinStoMap'}([m],[n])}\ar[rr]^{\ref{finalprop}}_{\sim}\ar[dd]^{F([m],[n])} && {\mathtt{FinStoMap'}([1],[n])^m}\ar[dd]^{\ref{Fsingleiso}}_\sim\\\\
{\mathtt{FinStoMap}([m],[n])}\ar[rr]^{\ref{matrixdecomp}}_{\sim} && {\mathtt{FinStoMap}([1],[n])^m}}
\eeq
which is commutative by construction of the maps. By proposition~\ref{Fsingleiso}, the right vertical arrow also is a bijection. Hence the diagram shows that the left vertical arrow also has to be bijective.
\end{proof}

\refs

\bibitem[Coxeter-Moser]{Cox} H. S. M. Coxeter, W. O. J. Moser, \textit{Generators and relations for discrete groups}, Springer (1957)
\bibitem[Fritz]{Fri} Tobias Fritz, \textit{Convex spaces I: definition and examples}, preprint. arXiv:0903.5522.
\bibitem[Lafont]{Laf} Yves Lafont, \textit{Towards an algebraic theory of Boolean circuits}, J. Pure Appl. Algebra 184, 257--310 (2003)
\bibitem[Massol]{Mas} A. Massol, \textit{Minimality of the system of seven equations for the category of finite sets}, Theoret. Comput. Sci. 176, 347--353 (1997)

\endrefs

\end{document}